\documentclass[11pt]{amsart}
\usepackage{latexsym,amssymb,amscd,amsthm,amsmath,amsfonts,exscale}
\usepackage[mathscr]{eucal}
\usepackage[all]{xy}
\usepackage{verbatim}

\def\dim{{\mathrm{dim}}}
\def\deg{{\mathrm{deg}}}
\def\Hzero{{\mathrm{H}^0}}

\def\max{{\mathrm{max}}}
\def\reg{{\mathrm{reg}}}
\def\sat{{\mathrm{sat}}}

\def\Proj{{\mathrm{Proj}}}

\def\Span{{\mathrm{Span}}}

\def\m{{\mathfrak m}}
\def\h{\mathrm{H}}

\def\KK{{\mathbb K}}
\def\PP{{\mathbb P}}

\def\III{{\mathscr I}}

\def\OOO{{\mathscr O}}

\def\EG{{\Xi}}

\def\ra{{\rightarrow}}

\theoremstyle{plain}
\newtheorem{ther}{Theorem}[section]
\newtheorem*{main}{Main Theorem}
\newtheorem{lem}[ther]{Lemma}
\newtheorem{prop}[ther]{Proposition}
\newtheorem{cor}[ther]{Corollary}

\theoremstyle{definition}
\newtheorem{defn}[ther]{Definition}

\theoremstyle{remark}

\begin{document}
\title{On the Castelnuovo-Mumford regularity of connected curves}
\author{Daniel Giaimo}
\address{Mathematics Department\\ UC Berkeley\\ Berkeley, California 94720}
\email{dgiaimo@math.berkeley.edu}
\date{30 August 2003}

\begin{abstract}
In this paper we prove the Eisenbud-Goto conjecture for connected
curves.  We also investigate the structure of connected curves for
which this bound is optimal.  In particular, we construct
connected curves of arbitrarily high degree in $\PP_\KK^4$ having
maximal regularity, but no extremal secants.  We also show that
any connected curve in $\PP_\KK^3$ of degree $\geq 5$ that has no
linear components and has maximal regularity has an extremal
secant.
\end{abstract}

\maketitle
    \setcounter{section}{-1}
    \section{Introduction} Let $S=\KK[x_0,...,x_n]$ where $\KK$ is an
        algebraically closed field, and let $M$ be a (graded)
        module over $S$.  We say that $M$ is $r$-regular if the
        $i$-th syzygy module of $M$ is generated in degrees less
        than or equal to $r+i$ for all $i$.  The regularity of $M$, denoted $\reg(M)$, is defined to be the
        infimum of all $r$'s such that $M$ is $r$-regular. If $X$ is a
        subscheme of $\PP_{\KK}^n$ then we define the regularity
        of $X$, denoted $\reg(X)$, to be the regularity of the (saturated homogeneous) ideal of
        $X$, $I_X$.

        In \cite{EG}, Eisenbud and Goto conjectured that if $X$ is
        a reduced nondegenerate connected in codimension 1 subscheme
        of $\PP_{\KK}^n$ then the regularity of $X$ should be no
        more than the degree of $X$ minus the codimension of $X$
        plus 1.  Various weakenings of this conjecture are known to be true in
        low dimension, (see, for example, \cite{GLP}, \cite{HP}, \cite{RL}, \cite{ZR}, \cite{SK1},
        \cite{SK2}, \cite{SK3}), but the main conjecture is still
        open.  In this paper we prove that
        this conjecture is true when $X$ is a curve.
        \begin{main}\label{T:main} Let $C$ be a connected reduced curve in $\PP_\KK^n$.
            Then
            \begin{equation}\label{T:EisenbudGoto'}
                \reg(C)\leq\deg(C)-\dim(\Span(C))+2
            \end{equation}
        \end{main}

        In \cite{GLP}, Gruson, Lazarsfeld, and Peskine proved this
        result in the case when $C$ is assumed to be irreducible.
        We prove our theorem by using their theorem as a base case
        and inducting on the number of irreducible
        components of $C$. The main lemma that allows this to work is
        a result of G. Caviglia \cite{GC}.

        We also study the connected curves for which the Eisenbud-Goto bound is optimal.
        In \cite{GLP}, Gruson, Lazarsfeld, and Peskine gave a
        complete classification of such curves in the irreducible
        case.  They proved that, except for a few low degree exceptions, all
        such curves are smooth, rational, and have a secant line of degree
        equal to the regularity of the curve.
        The existence of such a secant guarantees that the ideal
        of the curve will require a generator of degree equal to
        the regularity of the curve, so, at least for these curves, the
        regularity has to appear at the first step in a resolution
        for the ideal.

        In the connected case, it turns out that there are curves
        of maximal regularity which have
        arbitrarily high degree in $\PP_{\KK}^n$, $n\geq 4$ and which do not have
        such a secant.  A picture of such a curve of degree 9
        in $\PP_{\KK}^4$ follows.  (The labels of the various parts
        of the picture will be explained in section 3.)

$$
    
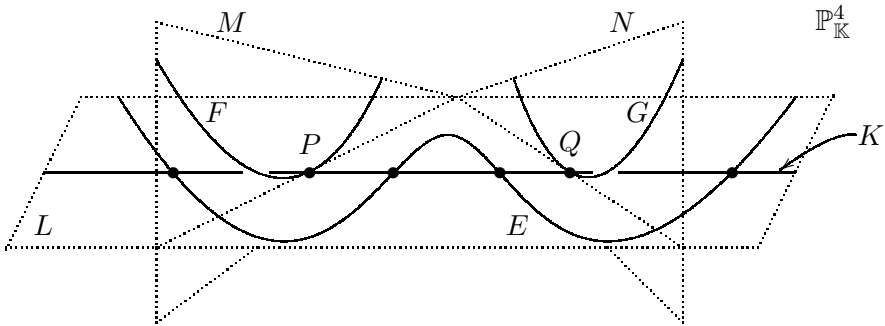
\xy
        (5,13)*{L},
        (115,25)*{K},
        (30,40)*{M},
        (68,13)*{E},
        (28,28)*{F},
        (84,28)*{G},
        (82,40)*{N},
        (110,40)*{\PP_\KK^4},
        (0,10)   **@{};
        (10,30)  **@{.};
        (110,30) **@{.};
        (100,10) **@{.};
        (0,10)   **@{.};
        (20,0)   **@{};
        (20,40)  **@{.};
        (60,30)  **@{.};
        (90,40)  **@{.};
        (90,0)   **@{.};
        (80,10)   **@{.};
        (33.33333,10)   **@{};
        (20, 0)  **@{.},
        (5, 20)  **@{};
        (31.5, 20) **@{-},
        (35, 20)  **@{};
        (78, 20) **@{-},
        (81.5, 20)  **@{};
        (105, 20) **@{-},
        (20,10)   **@{};
        (60,30)   **@{.};
        (90,10)   **@{.},
        (20,35);(50,32.5)**\crv{(37.4,4.8)},
        (90,35);(67.5,32.5)**\crv{(77,5)},
        (15,30);(105,30)**\crv{(40,-10)&(60,60)&(75,-10)},
        (40.4,23.8)*{P},
        (40.4,19.8)*{\bullet},
        (75,23.8)*{Q},
        (75,19.8)*{\bullet},
        (22.25,19.8)*{\bullet},
        (51.5,19.8)*{\bullet},
        (65.7,19.8)*{\bullet},
        (96.6,19.8)*{\bullet},
        (113,25);(103,20.4)**\crv{(110,25.5)} ?>*\dir{>}
    \endxy
$$

        The ideal of this curve, however, requires a generator of
        degree equal to the regularity, and we are still unsure as to whether
        the ideal of any curve of maximal regularity, except for a few low degree
        exceptions, requires a generator of degree equal to the regularity of the curve.

        In the final section we prove some structure theorems
        about connected curves of maximal regularity.  The first
        result that we prove here is that the irreducible components of such a
        curve have to be fairly well spaced apart.  To be precise,
        if $C$ is a connected curve that is the union of two connected curves $D$ and $E$
        having no common components then we are able to show that
        $\dim(\Span(D)\cap\Span(E))\leq 2$.  Using this fact, we
        split our analysis of these curves into three cases depending on whether
        $\dim(\Span(D)\cap\Span(E))$ equals 0, 1, or 2.  As a
        result of this analysis we are able to prove that any connected curve of maximal
        regularity in
        $\PP_\KK^3$ of degree $\geq 5$ having no
        linear components has an extremal
        secant.  We are also able to show that any connected
        subcurve of a connected curve of maximal regularity is either a line or
        has maximal regularity.

        \textbf{Acknowledgements.}  The author would like to thank
        his advisor, David Eisenbud, for suggesting he work on
        this problem, as well as for many useful conversations
        without which this paper would not be here today.

    \section{Notation}

    Throughout this paper we follow the following conventions:
        \begin{list}{$\bullet$}{}
            \item $\KK$ is an algebraically closed field of
            arbitrary characteristic.
            \item $S=\KK[x_0,...,x_n]$, $\m=(x_0,...,x_n)$, and $\PP_{\KK}^n=\Proj(S)$.
            \item All ideals and modules over $S$ are assumed to
            be graded.
            \item By a \emph{curve} we mean a reduced one-dimensional
            scheme embedded in $\PP_{\KK}^n$ for some $n$.
        \end{list}

    \section{Proof of the Main Theorem}
    To prove the main theorem, we take any reducible connected curve and write it
    as the union of two nonempty connected curves, $C$ and $D$, which meet each other in a finite set
    of points.  By induction on the number of irreducible components in the curves,
    we may assume that $C$ and $D$ satisfy the
    inequality (\ref{T:EisenbudGoto'}), so we need only show that this implies
    that $C\cup D$ does.  We will break the proof that $C\cup D$ satisfies (\ref{T:EisenbudGoto'}) into
    three cases depending on the dimension of $\Span(C)\cap\Span(D)$.  The
    following theorem of G. Caviglia allows us to bound the regularities
    of certain sums and intersections of ideals:
    \begin{ther}\label{T:Caviglia} Let $I$, $J\subset S$ be ideals
        and assume that $\dim(S/(I+J))\leq 1$.  Then
        $\reg(I+J)\leq\reg(I)+\reg(J)-1$, and $\reg(I\cap J)\leq\reg(I)+\reg(J)$.
    \end{ther}
    \begin{proof}
        See \cite{GC}.
    \end{proof}

    In order to simplify some of the statements that follow, we
    define $\EG(C)$ to be $\deg(C)-\dim(\Span(C))+2$ for any
    curve $C$.  With this definition the inequality (\ref{T:EisenbudGoto'})
    of the Main Theorem becomes
    \begin{equation}\label{T:EisenbudGoto}
        \reg(C)\leq\EG(C).
    \end{equation}
    The above theorem tells us, among other things, how the regularity of a
    union of two curves depends on the regularity of each of the
    component curves.  The following lemma, whose proof we
    leave to the reader, does the same
    thing for $\EG(C)$.
    \begin{lem}\label{T:EGBoundsum} Let $C$ and $D$ be curves with
        no common component.  Then $\EG(C\cup D)=\EG(C)+\EG(D)+
        \dim(\Span(C)\cap\Span(D))-2$ where we take the dimension
        of the empty scheme to be $-1$.
    \end{lem}

    Combining Theorem \ref{T:Caviglia} and Lemma \ref{T:EGBoundsum} we
    obtain the following corollary, which completes the proof
    of the inequality (\ref{T:EisenbudGoto}) for $C\cup D$
    assuming that $\dim(\Span(C)\cap\Span(D))\geq 2$.
    \begin{cor}\label{T:greaterthan2} Let $C$ and $D$ be curves with no common
        component, both satisfying (\ref{T:EisenbudGoto}).  Suppose
        that $\dim(\Span(C)\cap\Span(D))\geq 2$.  Then $C\cup D$
        satisfies (\ref{T:EisenbudGoto}).
    \end{cor}

    The second case we deal with is the case when $\dim(\Span(C)\cap\Span(D))=1$.
    To do this we need the following lemmas.
    \begin{lem}\label{T:intadd} Let $X$ and $Y$ be finite subschemes of
    $\PP_\KK^n$.  Then $\deg(X\cup Y)=\deg(X)+\deg(Y)-\deg(X\cap Y)$.
    \end{lem}
    \begin{proof}
        There is an exact sequence
        $$0\ra S_{X\cup Y}\ra S_X\oplus S_Y\ra S/(I_X+I_Y)\ra 0.$$
        The result follows as $S/(I_X+I_Y)$ is isomorphic to $S_{X\cap Y}$ in high
        degrees.
    \end{proof}
    \begin{lem}\label{T:min} If $C$ is a curve satisfying (\ref{T:EisenbudGoto})
    then $\EG(C)\geq 2$.
    \end{lem}
    \begin{proof}
        If $C$ is a line, then
        $$\EG(C)=\deg(C)-\dim(\Span(C))+2=1-1+2=2.$$
        If $C$ is not a line, then
        \[
        2\leq\reg(C)\leq\EG(C).\qedhere
        \]
    \end{proof}
    \begin{defn} We define the \emph{saturation degree}, denoted $\sat(I)$,
        of an ideal $I\subset S$
        to be $\reg(I^{sat}/I)+1$.  Equivalently it is the
        infimum of all $d$ for which $I_e=I^{sat}_e$ for all $e\geq d$.
    \end{defn}

    Notice that for any ideal $I\subset S$,
    \begin{equation}\label{regularity from saturation}
    \reg(I)=\max\{\reg(I^{sat}), \sat(I)\}.
    \end{equation}
    This follows from the long exact sequence in local cohomology coming
    from the short exact sequence
    $$0\ra I\ra I^{sat}\ra I^{sat}/I\ra 0$$
    after observing that
    $\h^0_\m(I^{sat})=\h^1_\m(I^{sat})=\h^i_\m(I^{sat}/I)=0$ for all $i\geq
    1$.
    Also notice that if $I$ and $J$ are ideals of $S$ then
    \begin{equation}\label{saturation of intersection}
        \sat(I\cap J)\leq\max\{\sat(I), \sat(J)\}
    \end{equation}
    as
    $(I\cap J)^{sat}=I^{sat}\cap J^{sat}$.

    We can now prove the following proposition which completes the proof
    of the inequality (\ref{T:EisenbudGoto}) for $C\cup D$
    assuming that $\dim(\Span(C)\cap\Span(D))=1$.
    \begin{prop}\label{T:equalto1} Let $C$ and $D$ be intersecting
        curves with no common component both satisfying
        (\ref{T:EisenbudGoto}).  Suppose that
        $\dim(\Span(C)\cap\Span(D))=1$.  Then $C\cup D$
        satisfies (\ref{T:EisenbudGoto}).
    \end{prop}
    \begin{proof}
        Let $L_1=\Span(C)$, $L_2=\Span(D)$, and
        $L=\Span(C)\cap\Span(D)$.  We separate the argument
        into two cases:

        \textbf{Case a:} $L$ is not a component of either $C$ or $D$.

        Let $J_1=I_C+I_L$ and $J_2=I_D+I_L$.  Notice that
        $$J_1+J_2=I_C+I_D+I_L=(I_C+I_{L_1})+(I_D+I_{L_2})=I_C+I_D$$
        as
        $I_{L_1}\subset I_C$ and $I_{L_2}\subset I_D$.  Therefore,
        from the exact sequence
        $$0\ra J_1\cap J_2\ra J_1\oplus J_2\ra J_1+J_2=I_C+I_D\ra 0$$
        we deduce that
        $$\reg(I_C+I_D)\leq\max\{\reg(J_1\cap J_2)-1,
        \reg(J_1),
        \reg(J_2)\}.$$
        By equation (\ref{regularity from saturation}) and inequality (\ref{saturation of intersection}),
        $$\reg(J_1\cap J_2)\leq\max\{\reg((J_1\cap
        J_2)^{sat}),\reg(J_1), \reg(J_2)\}.$$
        By Theorem
        \ref{T:Caviglia},
        \begin{align}
            \reg(J_1)&\leq\reg(C)\notag\\
            \reg(J_2)&\leq\reg(D)\notag.
        \end{align}
        Therefore,
        \begin{equation}\label{temp}
        \reg(I_C+I_D)\leq\max\{\reg((J_1\cap J_2)^{sat})-1, \reg(C),
        \reg(D)\}.
        \end{equation}

        The ideal of $(C\cap L)\cup
        (D\cap L)$ is $(J_1\cap J_2)^{sat}$.
        Since a finite scheme on a line has regularity
        equal to its degree,
        $$\reg((J_1\cap J_2)^{sat})=\deg((C\cap L)\cup (D\cap L)).$$
        By Lemma \ref{T:intadd},
        $$\deg((C\cap L)\cup (D\cap L))\leq\deg(C\cap L)+\deg(D\cap
        L)-\deg(C\cap D\cap L).$$
        Since $C$ and $D$ are assumed to
        intersect, and must do so in $L$,
        $$\deg(C\cap D\cap L)\geq 1.$$
        Since a scheme
        cannot have a secant of higher degree than its regularity,
        \begin{align}
            \deg(C\cap L)&\leq\reg(C)\notag\\
            \deg(D\cap L)&\leq\reg(D)\notag.
        \end{align}
        Therefore,
        $$\reg((J_1\cap
        J_2)^{sat})\leq\reg(C)+\reg(D)-1.$$
        Combining this with equation (\ref{temp}) it follows that
        \begin{equation}\label{temp2}
            \reg(I_C+I_D)\leq\reg(C)+\reg(D)-2.
        \end{equation}

        From equation (\ref{temp2}) and the exact sequence
        $$0\ra I_{C\cup D}\ra I_C\oplus I_D\ra I_C+I_D\ra 0$$
        it follows that
        $$\reg(C\cup D)\leq\reg(C)+\reg(D)-1.$$
        Because
        $C$ and $D$ both satisfy (\ref{T:EisenbudGoto}),
        $$\reg(C)+\reg(D)-1\leq\EG(C)+\EG(D)-1=\EG(C\cup
        D)$$
        by Lemma \ref{T:EGBoundsum}.  Therefore $C\cup D$ satisfies
        (\ref{T:EisenbudGoto}).

        \textbf{Case b:} $L$ is a component of one of $C$ or $D$.

        Without loss of generality, we assume that $L$ is a component of $C$.
        In this case $I_C\subset I_L$ and therefore
        $$I_C+I_D=I_L+I_D$$
        as $I_L\subset I_C+I_D$.
        Therefore
        $$\reg(I_C+I_D)=\reg(I_L+I_D)\leq\reg(D)$$
        by Theorem \ref{T:Caviglia}.
        From this and the exact sequence
        $$0\ra I_{C\cup D}\ra I_C\oplus I_D\ra I_C+I_D\ra 0$$
        we deduce that
        $$\reg(C\cup D)\leq\max\{\reg(C), \reg(D)+1\}.$$
        Because $C$ and $D$ both satisfy
        (\ref{T:EisenbudGoto}), this implies that
        $$\reg(C\cup D)\leq\max\{\EG(C), \EG(D)+1\}.$$
        We conclude that $C\cup D$ satisfies
        (\ref{T:EisenbudGoto}) by
        Lemma \ref{T:min} and Lemma \ref{T:EGBoundsum}.
    \end{proof}

    The final case we need to deal with to complete the proof of the
    main theorem is the case when $\dim(\Span(C)\cap\Span(D))=0$.
    (Notice that we cannot have $\dim(\Span(C)\cap\Span(D))=-1$ as
    $C$ and $D$ must intersect.)
    We do this in the following proposition.
    \begin{prop}\label{T:equalto0} Let $C$ and $D$ be intersecting
        curves with no common component both satisfying
        (\ref{T:EisenbudGoto}).  Suppose that
        $\dim(\Span(C)\cap\Span(D))=0$.  Then $C\cup D$
        satisfies (\ref{T:EisenbudGoto}).
    \end{prop}
    \begin{proof} Let $L_1=\Span(C)$, $L_2=\Span(D)$, and $P=\Span(C)\cap\Span(D)$.
    As in case a of Proposition \ref{T:equalto1},
    $$I_P=I_{L_1}+I_{L_2}\subset I_C+I_D\subset I_{C\cap D}.$$
    However, $C\cap D$ is nonempty, so its ideal cannot strictly contain a
    prime of dimension 1.  Therefore $I_P=I_C+I_D$.
    From the exact sequence
    $$0\ra I_{C\cup D}\ra I_C\oplus I_D\ra I_C+I_D=I_P\ra 0$$
    we see that
    $$\reg(C\cup D)\leq\max\{\reg(C), \reg(D),
    2\}.$$
    Because $C$ and $D$ both satisfy
    (\ref{T:EisenbudGoto}), this implies that
    $$\reg(C\cup D)\leq\max\{\EG(C), \EG(D), 2\}.$$
    We conclude that $C\cup D$ satisfies
    (\ref{T:EisenbudGoto}) by
    Lemma \ref{T:min} and Lemma \ref{T:EGBoundsum}.
    \end{proof}

    \section{Curves of Maximal Regularity: Introduction and Examples}

    In this section we partially classify
    those curves for which the Eisenbud-Goto bound is optimal.
    Such a curve is said to have \emph{maximal regularity}. The
    irreducible case is handled by Theorems 2.1 and 3.1 of \cite{GLP}.
    Their result is stated in terms of extremal secant lines which
    we define here.
   \begin{defn} A linear subspace $L\subset\PP_{\KK}^n$ is said to be an
        \emph{extremal secant} of a curve $C\subset\PP_{\KK}^n$
        if $C\cap L$ is finite and $\reg(C\cap L)=\EG(C)$.
    \end{defn}

    By Corollary 3.10 of \cite{EHU},
    if a connected curve $C$ has an extremal secant
    of dimension $i$ then $S/I_C$ has a minimal $i$-th
    syzygy of degree $\EG(C)+i$.  Hence, only curves of maximal regularity can have
    extremal secants.  What Gruson, Lazarsfeld, and Peskine showed is
    that essentially all irreducible curves of maximal regularity have extremal
    secant lines, and it follows from their classification that all
    irreducible curves of maximal regularity have extremal secant hyperplanes.

    We will prove that all connected curves which either have an extremal secant line or satisfy
    $\EG(C)=3$ have maximal regularity, but, unlike in the irreducible case, we will
    show how to construct connected curves of maximal regularity in $\PP_\KK^n$
    with no extremal secant line for any value of $\EG(C)\geq 4$ and
    any $n\geq 4$.  It would be interesting to know whether any such
    curves exist in $\PP_\KK^3$.  As we
    will see below, the construction we give necessitates
    $n\geq 4$.

    We have already shown that if a connected curve $C$ has an
    extremal secant line then it has maximal regularity.
    We prove in Proposition \ref{T:mincur} that
    all connected curves of regularity at most 2 have $\EG(C)=2$, so
    a connected curve with $\EG(C)=3$ must have regularity at least 3 hence exactly 3.
    Before proving this we need to define the
    following
    term which is used in the statement of the proposition.
    \begin{defn} We inductively define a \emph{linearly normal tree of rational curves}, shortened
    to \emph{tree} in the sequel.
    A linearly normal tree of rational curves with one component is defined to be a rational
    normal curve.  A linearly normal tree of rational curves with $k$ components, $k\geq 2$, is
    defined to be a curve $C$ in $\PP_{\KK}^n$ such that $C=C_0\cup D$ where $C_0$
    is a linearly normal tree of rational curves with $k-1$ components, $D$ is a rational
    normal curve and $C\cap
    D=\Span(C)\cap\Span(D)$ is a single point.
    \end{defn}
    \begin{prop}\label{T:mincur} Let $C$ be a connected curve in $\PP^n_{\KK}$.
                Then the following statements are equivalent:
        \begin{list}{}{}
            \item \textbf{(i)} $C$ is a tree.
            \item \textbf{(ii)} $\EG(C)=2$.
            \item \textbf{(iii)} $\reg(C)\leq 2$.
        \end{list}
    \end{prop}
    \begin{proof}  The equivalence of conditions \textbf{(i)} and
    \textbf{(ii)} follows directly from Xamb\'o's classification of
    connected in codimension 1 algebraic sets of minimal degree in
    \cite{SX}.  The fact that \textbf{(ii)} implies \textbf{(iii)}
    follows directly from the main theorem, so we need only prove
    that \textbf{(iii)} implies \textbf{(ii)}.

    If $\reg(C)=1$
    then $C$ is a line which satisfies
    $\EG(C)=2$, so assume that $\reg(C)=2$.
    Without loss of generality we may assume that $C$ is
    nondegenerate in $\PP^n_{\KK}$.  (Otherwise replace
    $\PP^n_{\KK}$ with $\Span(C)$.)
    From the exact sequence
    $$0\ra \III_C\ra \OOO_{\PP^n_{\KK}}\ra \OOO_C\ra 0$$
    we
    see that the map on global sections
    $\Gamma(\OOO_{\PP^n_{\KK}}(n))\ra\Gamma(\OOO_C(n))$ is
    surjective for all $n\geq 1$.  Since $C$ is connected it is
    surjective when $n=0$.  Therefore $C$ is projectively normal
    which implies $C$ is ACM as $C$ is a curve.  Therefore, if $Z=H\cap C$
    is a generic hyperplane section then
    $$\reg(I_{Z/H})=\reg(I_Z)=\reg(I_C)=2.$$
    From the exact sequence
    $$0\ra \III_{Z/H}\ra\OOO_H\ra\OOO_Z\ra 0$$
    we see that the map on global sections
    $$\KK^{n-1+1}\cong\Gamma(\OOO_H(1))\ra\Gamma(\OOO_Z(1))\cong
    \KK^{\deg(Z)}=\KK^{\deg(C)}$$
    is surjective.  Therefore
    $n\geq\deg(C)$.  Since $C$ is nondegenerate, $\deg(C)\geq n$ by Lemma \ref{T:min}.
    Therefore $\EG(C)=2$.
    \end{proof}
    \begin{cor} Let $C$ be a connected curve in $\PP_{\KK}^n$.  If
    $\EG(C)=3$ then $C$ is a curve of maximal regularity.
    \end{cor}

    We now construct examples of
    connected curves of maximal regularity with no extremal secant lines.
    By Proposition \ref{T:disjoint}, an
    example in $\PP_\KK^4$ gives rise to examples in
    $\PP_\KK^n$ for all $n\geq 4$.  Therefore we shall
    give the construction in $\PP_\KK^4$.  In fact, for the curves
    we construct in $\PP_\KK^4$ we will even show that they have
    no higher dimensional extremal secants as well.
    \begin{proof}[Construction]
    Let $m\geq 4$. Pick three 2-planes $L$, $M$, and $N$ such that $M$ and $N$ meet $L$
    in lines, but meet each other in a single point.  Pick a line $K$ in $L$
    meeting $M$ and $N$ in distinct points $P$ and $Q$.
    Pick a conic $F$ in $M$ which meets $L$ in a
    double point whose reduced structure is contained in $K$
    and a conic $G$ in $N$
    which meets $L$ in a double point whose reduced structure lies
    in $K$.  Pick a curve $E$ in $L$ of degree $m-3$ not meeting either
    $P$ or
    $Q$ and not containing $K$. Set $C=E\cup F\cup G\cup K$.
    (See the introduction for an example
    with $m=7$.)

    The span of $C$
    is all of $\PP_\KK^4$ so the regularity
    of $C$ is at most $\EG(C)=m-3+5-4+2=m$ as $\deg(C)=m-3+5$.
    We will show that $I_C$ requires a generator of degree at least $m$, hence has regularity
    at least $m$, and so exactly $m$.

    To prove this it suffices to show that any form $f\in I_C$ with
    $\deg(f)\leq m-1$ vanishes on the divisor $2K\subset L$.  Under
    the map $\KK[x_0, ...,
    x_4]\ra\KK[y_0, y_1, y_2]$ defining the embedding of $L$ in
    $\PP_\KK^4$, $f$ must map to an element, which we
    also call $f$, of $I_{C\cap
    L/L}$.  Let $g$ be the generator of $I_{E/L}$.  Then since
    $E\subset C\cap L$, we must have that $g|f$.  Let $h=f/g$. Let
    $D=F\cup G\cup K$.
    We can see that $h\in I_{D\cap L/L}$ as it still must vanish
    on $K$ and $f$ and $h$ differ only by a unit locally at $P$
    and $Q$.  We can also see that $h$ is either zero or has degree
    at most $2$.  In either case it must be a multiple of $k^2$
    where $k$ is the linear form defining $K$ as this is
    the only nonzero element of degree at most 2 in $I_{D\cap L/L}$ up
    to scalar multiples.  Therefore $f$ must vanish on $2K\subset L$.

    We must now prove that $C$ does not have an
    $m$-secant line.  To do this, let $H$ be any line in $\PP_{\KK}^4$ not contained in $C$.
    Suppose that $H$ does not lie in
    $L\cup M\cup N$.  Then it can meet $L\cup M\cup N$ in
    a scheme of length at most 3, hence it can be at most a 3-secant to
    $C$.  Since $m\geq 4$, such an $H$ is not an extremal secant to $C$.
    Suppose that $H$ lies in $M$ but not in $L$.
    Then it can meet $L\cup C$ in a scheme of length at most 3, hence it can
    be at most a 3-secant to $C$.  Again this implies that
    such an $H$ is not an extremal secant to $C$.
    By symmetry the same applies to any line which
    lies in $N$ but not $L$, so suppose
    that $H$ lies in $L$. Suppose that $H$ does not meet either $P$ or $Q$.  Then $H\cap C=H\cap
    (K\cup E)$ so $H$ is an $m-2$-secant to $C$. Again, such an $H$ is not an extremal secant to
    $C$, so suppose that $H$ passes through either $P$ or $Q$.  It
    cannot pass through both $P$ and $Q$ as then it would be equal
    to $K$, so it passes through exactly one of them which
    we may assume by symmetry is $P$.  Since $C\cap L$ has,
    locally at $P$, the structure of a line union a double point
    not in that line, $H\cap C$ has degree 2 locally
    at $P$.  Elsewhere, $H$ meets $C$
    only in $E$.  It must meet $E$ exactly $m-3$ times, so $H$ is an $m-1$
    secant to $C$.  Again, such an $H$ is not an extremal secant to
    $C$, so $C$ has no extremal secant lines.

    We now show that $C$ does not have an extremal
    secant 2-plane.  Suppose, to the contrary, that $H$ is such an
    extremal secant 2-plane.  By Lemma \ref{T:linearsection}, the
    degree of $C\cap H$ is at most $\EG(C)+1$.  Since
    $\EG(C)=m\geq 4$, Lemma \ref{T:finiteinplane} implies that
    there must be some line $M$ in $H$ such that $M\cap (C\cap H)$
    has degree $\EG(C)$.  But such an $M$ would be an extremal
    secant line to $C$ which we have already shown does not exist.
    Therefore $C$ does not have an extremal secant 2-plane.

    Finally we need to show that $C$ does not have an extremal
    secant hyperplane.  To do this we need the following lemma.
    \begin{lem}\label{construction lemma} Let $X$ be a finite subscheme of $\PP_\KK^3$ lying
    in a plane $O$ with $\deg(X)\geq 2$.  Let $Y$ be a finite subscheme of $\PP_\KK^3$
    of degree 2 which does not meet $O$.  Then $\reg(X\cup
    Y)\leq\reg(X)$.
    \end{lem}
    \begin{proof} The result follows from the exact sequence
    $$0\ra I_{X\cup Y}\ra I_X\oplus I_Y\ra I_X+I_Y\ra 0$$
    after observing that
    $I_X+I_Y$ is either $\m$ or $\langle x_0^2, x_1, x_2, x_3\rangle$ under a
    suitable choice of coordinates.
    \end{proof}

    Suppose that $C$ has an extremal secant hyperplane $H$.  If $H$
    passed through both $P$ and $Q$, it would
    contain $K$, which cannot happen as $H\cap C$ is finite.  By symmetry we may assume it does
    not contain $Q$.  In this case,
    $$H\cap C=(H\cap(E\cup F\cup K))\cup(H\cap G).$$
    Applying Lemma \ref{construction lemma} with $X=H\cap(E\cup F\cup
    K)$ and $Y=H\cap G$, we see that $H\cap\Span(E\cup F\cup K)$
    is an extremal secant 2-plane to $C$.  Since we have already
    shown that these don't exist, we have a contradiction.
    Therefore $C$ has no extremal secants.
    \end{proof}

    \section{Curves of Maximal Regularity: Structure Theorems}

    Despite the example above,
    it turns out that many
    connected curves of maximal regularity do have an
    extremal secant line.  One way one might try
    to show this is by breaking a connected curve into smaller
    connected components and seeing how they can fit together.
    The following corollary, which follows from Theorem \ref{T:Caviglia} and
    Lemma \ref{T:EGBoundsum}, says, in essence, that these
    components have to be fairly well separated.
    \begin{cor} Let $C$ and $D$ be intersecting connected curves
        in $\PP_{\KK}^n$ having no common components such that $C\cup D$
        has maximal regularity.  Then $\dim(\Span(C)\cap\Span(D))\leq
        2$.
    \end{cor}

    This corollary allows us to split our analysis of connected curves of
    maximal regularity into three major cases depending on
    $\dim(\Span(C)\cap\Span(D))$ where $C\cup D$ is the curve we
    are analyzing.  In case $\dim(\Span(C)\cap\Span(D))=0$, we have a complete
    classification which is the content of the following proposition.
    \begin{prop}\label{T:disjoint} Let $C$ and $D$ be intersecting connected
        curves in $\PP_{\KK}^n$ having no common components
        such that $\dim(\Span(C)\cap\Span(D))=0$,
        and assume, without loss of generality, that $\EG(C)\leq\EG(D)$.
        Then $C\cup D$ has maximal regularity if and only if $C$ is a tree
        of rational normal curves and $D$ is either a line or has maximal regularity.
        Moreover, $C\cup D$ has an extremal secant line if and only if
        $D$ is either a line or has an extremal secant line.
    \end{prop}
    \begin{proof}  Inspecting the proof of Proposition
    \ref{T:equalto0}, if $\EG(D)\geq 3$, then we see that there is
    equality throughout the chain of inequalities we proved if and only if
    $\reg(D)=\EG(D)$ and $\EG(C)=2$.  If $\EG(D)=2$, then
    $\EG(C)=2$ and
    $\EG(C\cup D)=2$, so by Proposition \ref{T:mincur} it follows that $C\cup D$
    has maximal regularity and $D$ is either a line or has maximal
    regularity. In either of these cases it follows from Proposition
    \ref{T:mincur} that $C$ is a tree.
    Therefore $C$ is a tree and $D$ is either a
    line or has maximal regularity.

    If $\EG(D)=2$, then, by the above, $C\cup D$ is a tree
    distinct from a line hence it has an
    extremal secant line.  Also, $D$ is either a line or has an
    extremal secant line.  Suppose that $\EG(D)\geq 3$.  If $D$ has
    an extremal secant line, then that same line is an extremal
    secant
    for $C\cup D$ as $\EG(D)=\EG(C\cup D)$.  Suppose that
    $C\cup D$ has an extremal secant line $L$.  If $L$ did not lie in either
    $\Span(C)$ or $\Span(D)$, then $L$ can intersect $C\cup D$ in at
    most a scheme of length 2, so $L$ must lie in one of $\Span(C)$
    or $\Span(D)$.  Notice that
    \begin{align}
        (C\cup D)&\cap\Span(C)=C\notag\\
        (C\cup D)&\cap\Span(D)=D\notag
    \end{align}
    so if $L$ lay in $\Span(C)$ it would
    intersect $C\cup D$ in a scheme of length at most 2.  Therefore $L$
    must lie in $\Span(D)$.  It follows that
    $$L\cap(C\cup D)=L\cap D$$
    and
    $L$ is an extremal secant for $D$.
    \end{proof}

    I like to think of this proposition as saying that
    one can add or remove ``feelers'' to or from a connected curve
    without changing whether it is a curve of maximal regularity
    where, by ``feelers'', I mean trees
    whose spans intersect the span of the rest of the curve in a
    single point.  Notice that this allows us to take the
    example constructed above and push it into any
    $\PP_\KK^n$, $n\geq 4$ by adding these
    ``feelers'' to it.

    Next suppose
    $\dim(\Span(C)\cap\Span(D))=1$.  This case splits naturally
    into two subcases depending on whether the line that is the
    intersection of the spans is contained in $C\cup D$ or not.
    We (partially) deal with the second case in the next
    proposition.
    \begin{prop}\label{T:dim1} Let $C$ and $D$ be intersecting connected curves in $\PP_{\KK}^n$
        having no common components such that
        $\dim(\Span(C)\cap\Span(D))=1$.  Let $L=\Span(C)\cap\Span(D)$, and assume
        that $L\nsubseteq C\cup D$.  If neither
        $C$ nor $D$ is a tree, then $C\cup D$ has maximal regularity
        if and only if $L$ is an extremal secant to $C\cup D$.
        Moreover, in this case $C$ and $D$ meet in a single point and
        $L$ is an extremal secant to
        both $C$ and $D$.
    \end{prop}
    \begin{proof} If $L$ is an extremal secant
        to $C\cup D$, then $C\cup D$ is a curve of maximal regularity.
        So we need only prove the opposite direction.

        Assume that $\reg(C\cup D)=\EG(C\cup D)$.  We use
        the notation from case a of Proposition
        \ref{T:equalto1}.  Looking at the end of the proof of that theorem
        we see that
        $$\reg(C\cup
        D)=\EG(C\cup D)=\EG(C)+\EG(D)-1$$
        implies that
        $$\max\{\EG(C), \EG(D),
        \reg(I_C+I_D)+1\}=\EG(C)+\EG(D)-1.$$
        Since both $\EG(C)$
        and $\EG(D)$ are greater than or equal to 2 this implies
        that
        $$\reg(I_C+I_D)=\EG(C)+\EG(D)-2.$$

        Looking further up the proof, we see that
        $$\reg(I_C+I_D)\leq\max\{\deg((C\cap L)\cup(D\cap L))-1, \EG(C),
        \EG(D)\}.$$
        Since neither $C$ nor $D$ is a tree, Proposition
        \ref{T:mincur} implies that $\EG(C)$ and $\EG(D)$ are both at
        least $3$, so
        \begin{equation}\label{temp9}
        \deg((C\cap L)\cup(D\cap L))=\EG(C)+\EG(D)-1=\EG(C\cup D).
        \end{equation}
        It follows that $L$ is an extremal secant to $C\cup D$ as
        $$(C\cap L)\cup(D\cap L)\subset (C\cup D)\cap
        L.$$
        Furthermore, since, by Lemma \ref{T:intadd},
        $$\deg((C\cap L)\cup(D\cap L))=\deg(C\cap L)+\deg(D\cap
        L)-\deg(C\cap D)$$
        it follows from equation (\ref{temp9}) that
        \begin{align}
            \deg(C\cap L)&=\EG(C)\notag\\
            \deg(D\cap L)&=\EG(D)\notag\\
            \deg(C\cap D)&=1.\notag
        \end{align}
        Therefore $C$ and $D$ meet in a single point and
        $L$ is an extremal secant to
        both $C$ and $D$.
    \end{proof}

    In one sense this proposition is more satisfactory than
    Proposition \ref{T:disjoint}, and in another it is less so.
    It is more satisfactory in the sense that if you have a curve
    that splits in this way, then you immediately deduce that the
    curve has an extremal secant without analyzing the subcurves
    any further.  It is less satisfactory because of the ``neither
    $C$ nor $D$ is a tree'' condition.  However,
    this condition cannot be removed.  If we
    break up the curve $C$ constructed at the end of the previous section
    as $F\cup(G\cup E\cup K)$ then
    the intersection of the spans is a line not contained in $C$
    even though $C$ doesn't have an extremal secant line.

    In fact, given any nonplanar connected curve $D$ with an
    extremal secant and any tree $C$, one can put these together
    in such a way that the intersection of their spans is a line
    which is not an extremal secant to $C\cup D$ even though $C\cup D$
    has extremal secant lines.  The following example
    demonstrates this.
    \begin{proof}[Construction] Let $N$ be a 3-plane in $\PP_\KK^4$ and
    let $D$ be a twisted cubic in $N$.  Let $L$ be a secant to $D$ and
    let $M$ be a line which intersects $D$ and $L$, but which is not
    a secant to $D$.  (This is actually true for every line other than $L$
    that passes through $D$ and $L$ and does so in distinct points.)   Let $C$ be a conic
    not contained in $N$ whose span contains $M$ and which intersects
    both $D$ and $L$.  A picture of such a configuration of curves follows.

$$
    
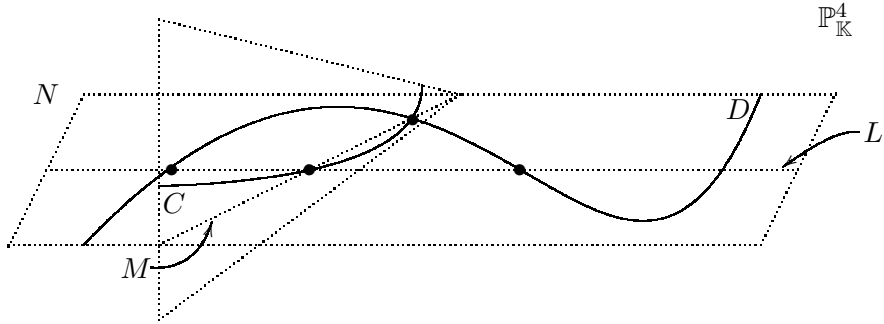
\xy
        (5,30)*=0{N},
        (17,7)*=0{M},
        (22,15.6)*=0{C},
        (97,28)*=0{D},
        (115,25)*=0{L},
        (110,40)*=0{\PP_\KK^4},
        (0,10)   **@{};
        (10,30)  **@{.};
        (110,30) **@{.};
        (100,10) **@{.};
        (0,10)   **@{.};
        (20,0)   **@{};
        (20,40)  **@{.};
        (60,30)  **@{.};
        (20, 0)  **@{.},
        (5, 20)  **@{};
        (105,20) **@{.},
        (20,10)  **@{};
        (60,30)  **@{.},
        (113,25);(103,20.4)**\crv{(110,25.5)} ?>*\dir{>},
        (19,7);(27,13)**\crv{(25,6.6)} ?>*\dir{>},
        (10,10);(100,30)**\crv{(60,62.7)&(80,-20)},
        (20,17.8);(55,31.3)**\crv{(55,19)},
        (21.7,19.85)*=0{\bullet},
        (40,19.85)*=0{\bullet},
        (68,19.85)*=0{\bullet},
        (53.7,26.6)*=0{\bullet},
     \endxy
$$

    $L$ is a 3-secant line to $C\cup D$ and
    $$\EG(C\cup D)=\EG(C)+\EG(D)-1=3$$
    so $L$ is an extremal secant to $C\cup
    D$.  However by applying Lemma \ref{T:linearsection} to the
    2-plane spanned by $L$ and $M$, which is the intersection of
    the spans of $C$ and $D$, we see that $M$ is only a 2-secant to
    $C\cup D$, hence not an extremal secant.
    \end{proof}

    The next proposition deals, (again, only partially), with the case when
    $\dim(\Span(C)\cap\Span(D))=1$ and $\Span(C)\cap\Span(D)\subset C\cup
    D$.
    \begin{prop} Let $C$ and $D$ be intersecting connected curves in $\PP_{\KK}^n$ having no
        common components such that
        $\dim(\Span(C)\cap\Span(D))=1$.  Let $L=\Span(C)\cap\Span(D)$, and assume
        that $L\subset C$.  Then $C\cup D$ has maximal regularity if and only if
        $C$ is a tree and $L\cup D$ has maximal regularity.
        Moreover, if $L\cup D$ has maximal regularity, then so
        does $D$.
    \end{prop}
    \begin{proof}
        First of all, suppose that $C$ is a tree and $L\cup D$
        has maximal regularity.  Then we can
        write $C=L\cup\bigcup_{i=1}^kC_i$ where the $C_i$'s are
        the connected components of $\overline{C\setminus L}$.  By
        Proposition \ref{T:mincur},
        all the $C_i$'s are trees and
        \begin{xalignat}{1}
        \dim(\Span(L\cup\bigcup_{i=1}^jC_i)\cap\Span(C_{j+1}))=0,\quad
        0\leq j\leq k-1.\notag
        \end{xalignat}
        Since $\Span(C)\cap\Span(D)=L$, it follows that
        \begin{xalignat}{1}
        \dim(\Span(D\cup
        L\cup\bigcup_{i=1}^jC_i)\cap\Span(C_{j+1}))=0, \quad 0\leq j\leq
        k-1.\notag
        \end{xalignat}
          By Proposition \ref{T:disjoint} and induction on $j$, it follows that $D\cup
        L\cup\bigcup_{i=1}^jC_i$ has maximal regularity for all $j$, $0\leq j\leq
        k$.  In particular, $C\cup D$ has maximal regularity.

        Now assume that $C\cup D$ has maximal regularity.  From
        the proof of case b of Proposition \ref{T:equalto1}, it follows that
        $$\max\{\EG(C), \EG(D)+1\}=\EG(C)+\EG(D)-1.$$
        Since $\EG(D)\geq 2$,
        this implies that $\EG(C)=2$.  By Proposition
        \ref{T:mincur} this implies that $C$ is a tree.

        To see that $L\cup D$ has maximal regularity, we let
        $C=L\cup\bigcup_{i=1}^kC_i$ where the $C_i$'s are
        the connected components of $\overline{C\setminus L}$.  As
        in the first part of this proof, all the
        $C_i$'s are trees and
        \begin{xalignat}{1}
        \dim(\Span(D\cup
        L\cup\bigcup_{i=1}^jC_i)\cap\Span(C_{j+1}))=0, \quad 0\leq j\leq
        k-1.\notag
        \end{xalignat}
        Therefore, by Proposition \ref{T:disjoint}
        and reverse induction on $j$, $D\cup
        L\cup\bigcup_{i=1}^jC_i$ has maximal regularity for all $j$, $0\leq j\leq
        k$.  In particular, $L\cup D$ has maximal regularity.

        Now suppose that $L\cup D$ has maximal regularity.  By
        Theorem \ref{T:Caviglia}, the regularity of $L\cup D$
        is at most 1 more than the regularity of $D$.  Since
        $\EG(L\cup D)=\EG(D)+1$, it follows that $D$
        must have maximal regularity.
    \end{proof}

    To finish the analysis we would need to answer the question,
    ``When can one add line to a curve of maximal
    regularity such that the resulting curve has maximal
    regularity?''  It might be tempting to guess that you have to add this line
    in such a way that it passes through an extremal secant line in a point that
    the original curve did not as
    this will always cause the resulting curve to have maximal
    regularity.  This is not the case as one can see from the
    example at the end of the last section by choosing $F$ to be the
    union of two lines $Y$ and $Z$ and splitting $C$ up as $Y\cup(Z\cup G\cup E\cup K)$.

    The last case is
    $\dim(\Span(C)\cap\Span(D))=2$.  We break this case into
    subcases depending on whether $C\cup D$ has components in
    $\Span(C)\cap\Span(D)$ or not.  We first deal with the case when no
    component of $C\cup D$ lies in $\Span(C)\cap\Span(D)$.  Before
    proceeding we need a few lemmas.
    \begin{lem} Let $C\subset\PP_{\KK}^n$ be a connected
        nondegenerate curve.  Let $H$ be a hyperplane containing no
        components of $C$.  Then $C\cap H$ is nondegenerate in $H$.
    \end{lem}
    \begin{proof} There is a commutative diagram with exact rows:

        \[
            \begin{CD}
                0 @>>> \Hzero(\OOO_{\PP_\KK^n})  @>>> \Hzero(\OOO_{\PP_\KK^n}(1)) @>>> \Hzero(\OOO_H(1)) @>>> 0\\
                @.     @V\alpha VV                   @V\beta VV                   @V\gamma VV\\
                0 @>>> \Hzero(\OOO_C)            @>>> \Hzero(\OOO_C(1))           @>>> \Hzero(\OOO_{C\cap H}(1))\\
            \end{CD}
        \]

        $C$ being nondegenerate is equivalent to the statement
        that $\beta$ is injective.  $C$ being connected is
        equivalent to the statement that $\alpha$ is surjective.
        It follows from the snake lemma that $\gamma$ is
        injective.  This says precisely that $C\cap H$ is
        nondegenerate in $H$.
    \end{proof}
    \begin{lem}\label{T:linearsection} Let $C\subset\PP_{\KK}^n$ be a connected
        nondegenerate curve.  Let $L$ be a linear subspace containing no
        components of $C$.  Then $$\deg(C\cap
        L)\leq\deg(C)-n+1+\dim(L)=\EG(C)+\dim(L)-1.$$
    \end{lem}
    \begin{proof} Choose a hyperplane $H$ containing $L$ but not
    containing any component of $C$, and consider the
    following commutative triangle.
    $$\xymatrix{\Hzero(\OOO_H(1))\ar[r]^{\alpha}\ar[dr]^{\beta}
              & \Hzero(\OOO_{C\cap H}(1))\ar[d]^{\gamma} \\ &
                \Hzero(\OOO_{C\cap L}(1))}$$
    Since $C\cap H$ is nondegenerate in $H$, the map $\alpha$ is
    injective.  Since $C\cap L\subset L$,
    $$\dim_\KK(\ker(\beta))\geq\dim(H)-\dim(L).$$
    Therefore
    $$\dim_\KK(\ker(\gamma))\geq\dim(H)-\dim(L).$$
    Since $C\cap H$ is
    a finite scheme, $\gamma$ is surjective.  Therefore
    \begin{equation}
    \begin{split}
    \deg(C\cap L)&=\dim_\KK(\Hzero(\OOO_{C\cap L}(1)))\\
    &=\dim_\KK(\Hzero(\OOO_{C\cap H}(1)))-\dim_\KK(\ker(\gamma))\\
    &\leq\deg(C)-(n-1)+\dim(L).\notag
    \end{split}
    \end{equation}
    \end{proof}
    \begin{lem}\label{T:finiteinplane} Let $X\subset\PP_{\KK}^2$ be a finite scheme of
        degree $d$.  Then $\reg(X)\leq d$.  Moreover, if $\reg(X)=d$
        then $X$ lies on a line, $M\subset\PP_{\KK}^2$. If
        $\reg(X)=d-1$ and $d\neq 4$, then there exists a line $M\subset\PP_{\KK}^2$
        such that $\deg(M\cap X)=d-1$.
    \end{lem}
    \begin{proof}  By Proposition 3.7 and Corollary 3.9 of \cite{DE}, if
        $$0\ra\sum_{i=1}^tS(-b_i)\xrightarrow{M}\sum_{i=1}^{t+1}S(-a_i)\ra
        S$$
        is a free resolution of $S_X$ with $a_1\geq ...\geq
        a_{t+1}$ and $b_1\geq ... \geq b_t$, and if we let $e_i$
        and $f_i$ denote the degrees of the entries on the
        principal diagonals of $M$, then for all $i$,

        $e_i\geq 1$, $f_i\geq 1$,

        $f_i\geq e_i$, $f_i\geq e_{i+1}$,

        $a_i=\sum_{j<i}e_j+\sum_{j\geq i}f_j$,

        $b_i=a_i+e_i$, for $1\leq i\leq t$ and
        $\sum_{i=1}^tb_i=\sum_{i=1}^{t+1}a_i$,

        $d=\sum_{i\leq j}e_if_j$.

        Notice that this implies that
        $$b_1-1=a_1+e_1-1\geq a_1.$$
        Therefore,
        $$\reg(X)=\max\{a_1,
        b_1-1\}=b_1-1=a_1+e_1-1=\sum_{i=1}^tf_t+e_1-1.$$
        Suppose
        that $t\geq 2$.  Then
        $$\sum_{i=1}^tf_t+e_1-1\leq\sum_{i=1}^te_1f_t<\sum_{i\leq
        j}e_if_j=d.$$
        Suppose that $t=1$. Then
        $$\sum_{i=1}^tf_t+e_1-1=f_1+e_1-1\leq e_1f_1=d.$$
        Therefore $\reg(X)\leq d$.

        Suppose that $\reg(X)=d$.  By the above argument,
        $t=1$ and one of $e_1$ or $f_1$ equal to 1.
        Since $f_1\geq e_1$ this implies that $e_1=1$.  Therefore
        $a_2=1$ and there is some line
        $M\subset\PP_{\KK}^2$ such that $X\subset M$.

        Suppose that $\reg(X)=d-1$ and $d\neq 4$.  If $d=3$ we are done as
        $X$ has a degree 2 subscheme which must then lie on a line.
        Therefore we may assume that $d\geq 5$. By the
        above argument, either $t=1$ and
        $e_1=f_1=2$, or $t\geq 2$ and
        \begin{equation}\label{temp3}
        \sum_{i=1}^tf_t+e_1=\sum_{i=1}^te_1f_t+1=\sum_{i\leq
        j}e_if_j.
        \end{equation}
        In the first case we would have $d=e_1f_1=4$
        which we already ruled out by assuming $d\geq 5$.
        Therefore we must be in the second case.
        In this case, equation (\ref{temp3}) implies that
        \begin{equation}
        e_1=1\qquad
        t=2  \qquad
        e_2=f_2=1.\notag
        \end{equation}
        Also, we must have that $b_1=d$ since
        $\reg(X)=b_1-1$.  Therefore
        \begin{equation}
        a_1=b_1-e_1=d-1\qquad a_2=e_1+f_2=2.\notag
        \end{equation}
        $X$ cannot lie on a line as then it
        would have regularity $d$.  Therefore $a_3\geq 2$, so
        $a_3=2$ as $2=a_2\geq a_3$.  Also,
        $$d+b_2=b_1+b_2=a_1+a_2+a_3=d-1+2+2=d+3$$
        so $b_2=3$.
        Therefore the free resolution of $S_X$ has the form:
        $$0\ra S(-d)\oplus S(-3)\ra S(-d+1)\oplus S(-2)\oplus S(-2)\ra S.$$

        Let $f$, $g$, and $h$ form a generating set for $I_X$ with
        $f$ and $g$ of degree 2 and $h$ of degree $d-1$.
        If $f$ is a nonzerodivisor modulo $g$, then $Z((f,g))$
        would be a complete intersection of degree 4 containing
        $X$.  This is a contradiction to the assumption that
        $\deg(X)\geq 5$, so $f$ must be a zerodivisor modulo $g$.
        The only way this can happen is if $f=ml_1$ and $g=ml_2$
        for some linear forms $m$, $l_1$, and $l_2$ with $l_1\neq
        l_2$. If $h\in (m)$ then
        the line $M$ defined by $m$ would lie inside of $X$.
        Therefore $h$ is a nonzerodivisor modulo $m$ and
        $Z((h,m))=M\cap X$ is a finite scheme of degree $d-1$.
    \end{proof}

    We return to the case
    $\dim(\Span(C)\cap\Span(D))=2$ and no component of $C\cup D$
    lies in $\Span(C)\cap\Span(D)$.  Like Proposition
    \ref{T:dim1}, it turns out that
    such a curve of maximal regularity must have an extremal
    secant without looking any further at subcurves.  In this
    case we don't have the restriction that
    neither of the components be trees.
    \begin{prop} \label{T:dim2}Let $C$ and $D$ be intersecting connected curves in $\PP_{\KK}^n$
        having no common components such that
        $L=\Span(C)\cap\Span(D)$ has dimension $2$, and assume that no component of
        either $C$ or $D$ lies in $L$.  Then $C\cup D$ has maximal regularity if and only if
        there is some line $M\subset L$ which is an extremal
        secant for $C\cup D$.
    \end{prop}
    \begin{proof}
        If $C\cup D$ has an extremal
        secant lying in $L$ then it has maximal regularity, so we
        need only prove the reverse direction.

        Assume that $\reg(C\cup D)=\EG(C\cup D)$.  Then from the
        exact sequence
        $$0\ra I_{C\cup D}\ra I_C\oplus I_D\ra I_C+I_D\ra 0$$
        we see that
        $$\reg(C\cup
        D)\leq\max\{\reg(C), \reg(D), \reg(I_C+I_D)+1\}.$$
        By the Main Theorem we have that $\reg(C)\leq\EG(C)$ and
        $\reg(D)\leq\EG(D)$, so from
        Lemma \ref{T:EGBoundsum} we see that
        $$\reg(I_C+I_D)\geq\EG(C\cup D)-1.$$
        By Theorem \ref{T:Caviglia}
        $$\reg(I_C+I_D)\leq\reg(C)+\reg(D)-1\leq\EG(C\cup D)-1.$$
        Therefore
        \begin{equation}\label{temp4}
        \reg(I_C+I_D)=\EG(C\cup D)-1.
        \end{equation}

        Now let $J_1=I_C+I_L$ and $J_2=I_D+I_L$.  Since
        $I_L\subset I_C+I_D$, we see that $J_1+J_2=I_C+I_D$.
        From the exact sequence
        $$0\ra J_1\cap J_2\ra J_1\oplus J_2\ra J_1+J_2=I_C+I_D\ra 0$$
        it follows that
        $$\reg(I_C+I_D)\leq\max\{\reg(J_1), \reg(J_2),
        \reg(J_1\cap J_2)-1\}.$$
        By Theorem \ref{T:Caviglia} and the Main Theorem
        \begin{align}\label{temp5}
            \reg(J_1)&\leq\reg(C)\leq\EG(C)\\\label{temp6}
            \reg(J_2)&\leq\reg(D)\leq\EG(D).
        \end{align}
        Therefore
        $$\reg(I_C+I_D)\leq\max\{\EG(C), \EG(D),\reg(J_1\cap
        J_2)-1\}.$$
        By equation (\ref{temp4}), this implies that
        $$\reg(J_1\cap J_2)\geq\EG(C\cup D).$$
        By
        Theorem \ref{T:Caviglia} this inequality can't be
        strict, so
        \begin{equation}\label{temp7}
        \reg(J_1\cap J_2)=\EG(C\cup D).
        \end{equation}

        By equation (\ref{regularity from saturation}),
        $$\reg(J_1\cap J_2)=\max\{\reg((J_1\cap J_2)^{sat}),
        \sat(J_1\cap J_2)\}.$$
        From this, inequalities (\ref{saturation of
        intersection}), (\ref{temp5}), and (\ref{temp6}), and equation (\ref{temp7}) we see
        that
        $$\reg((J_1\cap
        J_2)^{sat})=\EG(C\cup D).$$
        The ideal of $(C\cap L)\cup
        (D\cap L)$ is $(J_1\cap J_2)^{sat}$.
        Since $(C\cap
        L)\cup (D\cap L)$ lies in $(C\cup D)\cap L$,
        \begin{equation}
        \deg((C\cap L)\cup (D\cap L))\leq\deg((C\cup D)\cap
        L).\notag
        \end{equation}
        By Lemma
        \ref{T:linearsection}
        \begin{equation}
        \deg((C\cup D)\cap L)\leq\EG(C\cup D)+1.\notag
        \end{equation}
        Therefore $(C\cap
        L)\cup (D\cap L)$ is a finite subscheme of a 2-plane, $L$,
        of degree at most $\EG(C\cup D)+1$ and regularity exactly
        $\EG(C\cup D)$.  Furthermore,
        $$\EG(C\cup
        D)=\EG(C)+\EG(D)\geq 4$$
        so either the degree of $(C\cap
        L)\cup (D\cap L)$ is exactly equal to its regularity, or
        it is at least 5 and the difference between them is 1.  In either case, Lemma
        \ref{T:finiteinplane} implies that there is a line
        $M\subset L$ such that
        $$\deg(M\cap((C\cap
        L)\cup (D\cap L)))=\EG(C\cup D).$$
        This implies
        that the length of $M\cap(C\cup D)$
        is at least $\EG(C\cup D)$, hence $M$ is an extremal secant line to $C\cup D$.
    \end{proof}

    We now turn our attention to the case when
    $\dim(\Span(C)\cap\Span(D))=2$ and there are components of $C\cup D$
    that lie in $\Span(C)\cap\Span(D)$.
    We need two lemmas
    about subschemes of $\PP_\KK^2$.
    \begin{lem} \label{T:planarregularity}Let $X=D\cup Y$ be a subscheme of $\PP_\KK^2$
    where $D$ is a curve of degree $d$ and $Y$ is a finite scheme.  If
    $$0\ra\sum_{i=1}^tS(-b_i)\xrightarrow{M}\sum_{i=1}^{t+1}S(-a_i)\ra
    S$$
    is a free resolution of $S_X$ with $a_1\geq ...\geq a_{t+1}$ and
    $b_1\geq ... \geq b_t$, and if we let $e_i$ and $f_i$ denote the
    degrees of the entries on the principal diagonals of $M$, then for all $i$,

    $e_i\geq 1$, $f_i\geq 1$,

    $f_i\geq e_i$, $f_i\geq e_{i+1}$,

    $a_i=\sum_{j<i}e_j+\sum_{j\geq i}f_j+d$,

    $b_i=a_i+e_i$, for $1\leq i\leq t$ and
    $\sum_{i=1}^tb_i+d=\sum_{i=1}^{t+1}a_i$,

    $H_{S_X}(n)=dn+1-(d-1)(d-2)/2+\sum_{i\leq j}e_if_j$.
    \end{lem}

    The proof of this lemma is analogous to Proposition 3.7
    and Corollary 3.9 of \cite{DE} and we do not repeat it here.
    The main difference is that the generators of the ideal
    are not the maximal minors of $M$ but rather the maximal
    minors of $M$ multiplied by the equation of $D$ which is where
    you get all the ``$+d$'' terms above as well as the difference in
    the Hilbert polynomial.
    \begin{lem}\label{T:notsofiniteinplane} Let $\emptyset\neq
        D\subset\PP_{\KK}^2$ be a curve of
        degree $d$ and let $Y\subset\PP_{\KK}^2$ be a finite scheme.
        Then
        $$\reg(D\cup Y)\leq d+\deg(Y)-\deg(D\cap Y).$$
        Moreover,
        if we have equality and $\deg(D\cap Y)=1$, then there is
        some line, $M\subset\PP_{\KK}^2$, such that
        $$\deg(M\cap (D\cup Y))=d+\deg(Y)-1.$$
    \end{lem}
    \begin{proof}  If $Y\subset D$ then everything is trivial, so
    we may assume that $Y\nsubseteq D$.
    From the
    exact sequence
    $$0\ra S_{D\cup Y}\ra S_D\oplus S_Y\ra S/(I_D+I_Y)\ra 0$$
    and the fact that $S/(I_D+I_Y)$
    agrees with $S_{D\cap Y}$ in high degrees we see that
    \begin{equation}
    \begin{split}
    H_{S_{D\cup
    Y}}(n)&=H_{S_D}(n)+H_{S_Y}(n)-H_{S_{D\cap Y}}(n)\\
    &=dn+1-(d-1)(d-2)/2+\deg(Y)-\deg({D\cap Y}).\notag
    \end{split}
    \end{equation}
    Therefore, using the notation
    from Lemma \ref{T:planarregularity} applied to $D\cup Y$,
    $$\sum_{i\leq j}e_if_j=\deg(Y)-\deg({D\cap
        Y}).$$
        Also,
    \begin{equation}
    \begin{split}
        \reg(D\cup Y)&=\max\{a_1,
        b_1-1\}=b_1-1\\&=a_1+e_1-1=\sum_{i=1}^tf_t+e_1+d-1.\notag
    \end{split}
    \end{equation}
        If $t\geq 2$, then
        $$\sum_{i=1}^tf_t+e_1+d-1\leq\sum_{i=1}^te_1f_t+d<\sum_{i\leq
        j}e_if_j+d=d+\deg(Y)-\deg({D\cap
        Y}).$$
        If $t=1$, then
        $$\sum_{i=1}^tf_t+e_1+d-1=f_1+e_1+d-1\leq e_1f_1+d=d+\deg(Y)-\deg({D\cap
        Y}).$$
        Therefore
        $\reg(D\cup Y)\leq d+\deg(Y)-\deg({D\cap
        Y})$.

        Suppose that $\deg(D\cap Y)=1$ and $\reg(D\cup Y)=d+\deg(Y)-1$.
        By the above argument, this implies that $t=1$ and one of $e_1$ or $f_1$ equal to 1.
        Since $f_1\geq e_1$ this implies that $e_1=1$.  Therefore
        $a_2=d+e_1=d+1$.  Also, we must have that $b_1=d+\deg(Y)$, so the
        free resolution of $S_{D\cup Y}$ is
        \begin{equation}
        0\ra S(-(d+\deg(Y)))\ra S(-(d+\deg(Y))+1)\oplus S(-(d+1))\ra
        S\notag
        \end{equation}
        Let $g$ be the equation of $D$.  Then the above says that
        the minimal generators of $I_{D\cup Y}$ have the form $gm$ and
        $gh$ where $m$ is some linear form and $h$ is a form of
        degree $\deg(Y)-1$.  Moreover, $m\nmid h$ as otherwise
        $gh$ wouldn't be a minimal generator.  Therefore, if
        $m\nmid g$, then
        $Z((gh,m))=M\cap (D\cup Y)$ is a finite scheme of degree
        $d+\deg(Y)-1$ where $M$ denotes the line defined by $m$.

        Suppose that $m|g$.  Since $\deg(D\cap
        Y)=1$ and $Y\nsubseteq D$, it follows that $Y$ must be a
        double point not lying in $D$ whose reduced point is
        a nonsingular point of $D$ lying on $M$. In this case, $h$ is a linear
        form which does not divide $g$.  Therefore, as above, $h$ defines a
        line which intersects $D\cup Y$ in a finite scheme of
        length $d+1$.
    \end{proof}
    \begin{prop} \label{T:planarD}Let $C$ and $D$ be intersecting connected curves in
        $\PP_{\KK}^n$ having no common components
        such that $L=\Span(C)\cap\Span(D)$ has dimension $2$. Assume that no
        component of $C$ lies in $L$ and that $D\subset L$.  Then $C\cup D$
        has maximal regularity if and only if
        there is some line $M\subset L$ which is an extremal
        secant for $C\cup D$.
    \end{prop}
    \begin{proof}
        If $C\cup D$ has an extremal
        secant lying in $L$ then it has maximal regularity, so we
        need only prove the reverse direction.

        Assume that $\reg(C\cup D)=\EG(C\cup D)$.  Then exactly as
        in the proof of Proposition \ref{T:dim2}, we have that
        $$\reg(I_C+I_D)=\EG(C\cup D)-1.$$
        Let $J=I_C+I_L$.
        Since
        $I_L\subset I_C+I_D$, we see that $J+I_D=I_C+I_D$.
        Again, as in the proof of Proposition \ref{T:dim2},
        from the exact sequence
        $$0\ra J\cap I_D\ra J\oplus I_D\ra J+I_D=I_C+I_D\ra 0$$
        we
        deduce that
        $$\reg((J\cap
        I_D)^{sat})=\EG(C\cup D).$$
        Since $Z((J\cap I_D)^{sat})=D\cup (C\cap L)$, Lemma \ref{T:notsofiniteinplane}
        implies that
        \begin{equation}\label{temp8}
        \deg(C\cap D)\leq \deg(D)+\deg(C\cap L)-\EG(C\cup D).
        \end{equation}
        By Lemma \ref{T:linearsection},
        \begin{equation}\label{temp12}
        \deg(C\cap L)\leq\EG(C)+1.
        \end{equation}
        Also,
        \begin{align}
        \EG(D)&=\deg(D)\notag\\
        \EG(C\cup
        D)&=\EG(C)+\EG(D).\notag
        \end{align}
        Therefore $\deg(C\cap D)=1$.
        Moreover, we must have equality in (\ref{temp8}) and (\ref{temp12}) and thus by Lemma \ref{T:notsofiniteinplane},
        there is a line $M\subset L$, such that
        \begin{equation}
        \deg(M\cap(D\cup(C\cap L)))=\EG(C\cup D).\notag
        \end{equation}
        Since
        $$D\cup(C\cap L)\subset(C\cup D)\cap L$$
        $M$ is an extremal secant to $C\cup D$.
    \end{proof}

    These results allow us to prove that many
    of the curves of maximal regularity in $\PP_\KK^3$ have extremal secants.
    \begin{ther} Let $C\subset\PP_\KK^3$ be a connected curve with
        no linear components.  Then $C$ has maximal regularity if
        and only if either $\EG(C)=3$ or $C$ has an extremal
        secant.
    \end{ther}
    \begin{proof} We have already proved the ``if'' part of this
    theorem, so assume that $C$ has maximal regularity.  If $C$ is
    planar then $C$ has an extremal secant.

    Assume that
    $C$ is nondegenerate in $\PP_\KK^3$.  If $C$ is irreducible,
    then we are done by Theorems 2.1 and 3.1 of \cite{GLP}.

    Assume
    that $C$ is reducible.  We can find two connected
    curves with no common components, $D$ and $E$, such that $C=D\cup E$.  In fact, if $A$
    and $B$ are planar components of $C$ lying in the same plane
    then we can (and do) ensure that $A$ and $B$ either both lie in $D$ or
    both lie in $E$.
    If $D$
    is nondegenerate in $\PP_\KK^3$ and $E$ is planar then
    Proposition \ref{T:planarD} implies that $C$ has an
    extremal secant in the plane that $E$ spans.

    Assume that $D$ and $E$ are both
    planar and the intersection of their spans, $L$, is a line lying in
    neither of them.  If neither $D$ nor $E$ is a conic then
    Proposition \ref{T:dim1} implies that $L$ is an extremal
    secant to $C$.

    Assume that $D$ is a conic.  If $\deg(D\cap
    E)=1$, then
    $$\deg((D\cap L)\cup(E\cap
    L))=\deg(D)+\deg(E)-1=\EG(C).$$  Since
    $$(D\cap L)\cup(E\cap
    L)\subset C\cap L$$
    this implies that $L$ is an extremal
    secant to $C$.

    Assume that
    $\deg(D\cap E)\geq 2$.  Since $D$ is a conic, $\deg(D\cap
    L)=2$, so $\deg(D\cap E)=2$.  Notice that $I_D+I_E$
    contains $I_L$.  It also contains a polynomial, $f$, of degree 2
    which defines $D$ in its span.  However, $I_L+(f)$ is the
    (saturated) homogeneous ideal of a scheme of length 2 on $L$.
    Therefore, since $D\cap E$ is a subscheme of this scheme and
    also has length 2, we must have equality and $I_D+I_E=I_{D\cap
    E}$.  In particular, since $D$, $E$, and $D\cap E$ are all
    ACM, $C$ is ACM.  We
    now let $H$ be any plane such that $Z=C\cap H$ is finite.
    Since $C$ is ACM,
    $\reg(Z)=\reg(C)$.  If $\deg(C)=4$, then $\EG(C)=3$ and we are
    done.

    Assume that $\deg(C)\geq 5$.  Then
    $\deg(Z)\geq 5$, and Proposition \ref{T:finiteinplane} implies
    that there is some line $M$ such that $\deg(M\cap Z)=\EG(C)$.
    In particular, $\deg(M\cap C)=\EG(C)$ and $M$ is an extremal
    secant to $C$.  This finishes the proof.
    \end{proof}

    Finally, we mention the following corollaries of
    the analysis in this section.
    \begin{cor} Let $C$ and $D$ be intersecting connected curves in
        $\PP_{\KK}^n$ having no common components, and suppose that $C\cup D$ has maximal
        regularity.  Then $C$ and $D$ are both either lines or curves having maximal
        regularity.
    \end{cor}
    \begin{cor}  Any connected subcurve of a curve of maximal
    regularity is either a line or has maximal regularity.
    \end{cor}

\bibliographystyle{amsplain}
\bibliography{CC}
\end{document}